\documentclass[psamsfonts, reqno]{amsart}
\usepackage{amsfonts,amssymb,amsmath,amscd,amsthm}
\usepackage[all]{xy}

\providecommand{\cprod}[2]{\ensuremath{ #1\rtimes_{#2}\mathbb{Z}}}
\providecommand{\rcprod}[2]{\ensuremath{#1\rtimes_{#2}\mathbb{R}}}

\providecommand{\ts}[1]{{\textstyle{#1}}}

\newcommand{\C}{{\mathbb C}}

\newcommand{\Z}{{\mathbb Z}}
\newcommand{\R}{\mathbb R}
\newcommand{\Q}{{\mathbb Q}}
\newcommand{\T}{{\mathbb T}}

\newcommand{\Aut}{\mathrm{Aut}}

\newcommand{\calA}{\mathcal{A}}

\newcommand{{\M}}{\mathcal{M}}
\newcommand{\ttheta}{\widetilde{\theta}}
\newcommand{\talpha}{\widetilde{\alpha}}
\newcommand{{\frakS}}{\mathfrak{S}}

\renewcommand{\phi}{\varphi}
\newcommand{\E}{{\Lambda}}

\newcommand{\CSpec}{\Gamma}

\newtheorem{theorem}{Theorem}[section]
\newtheorem{lemma}[theorem]{Lemma}
\newtheorem{corollary}[theorem]{Corollary}
\newtheorem{proposition}[theorem]{Proposition}

\begin{document}

\title{Continuous and discrete flows on operator algebras}
\author{Benjam\'{\i}n Itz\'a-Ortiz}
%\date\today

\begin{abstract}
Let $(N,\R,\theta)$ be a centrally ergodic W* dynamical system. When $N$ is not a factor, we show that, for each $t\not=0$,  the crossed product induced by the time $t$ automorphism $\theta_t$ is not a factor if and only if there exist a rational number $r$ and an eigenvalue $s$ of the restriction of $\theta$ to the center of $N$, such that $rst=2\pi$. In the C* setting, minimality seems to be  the notion corresponding to central ergodicity. We show that if $(A,\R,\alpha)$ is a minimal unital C* dynamical system and $A$ is either prime or commutative but not simple, then, for each $t\not=0$, the crossed product induced by the time $t$  automorphism $\alpha_t$ is not simple if and only if there exist a rational number $r$ and an eigenvalue $s$ of the restriction of $\alpha$ to the center of $A$, such that $rst=2\pi$.
\end{abstract}

\maketitle

\section*{Introduction}
Recall that a flow $(Y,T)$ is a pair consisting of a compact metric space $Y$ and an action $T\colon Y\times \R\to Y$. The time $t$ map of the flow $(Y,T)$ is the automorphism $T^t\colon Y\to Y$. We say that a flow $(Y,T)$ is minimal if there is no nontrivial closed invariant subspaces of $Y$.

If $(Y,T)$ is a minimal flow then, for $t\not=0$, Proposition~1.5 in  \cite{IO} shows that the crossed product induced by the time $1/t$ map $T^{\frac{1}{t}}$ is not minimal if and only if there exists a rational number $r$ and an eigenvalue $s$ of $T$ such that $rst=2\pi$. (The $2\pi$ term appears in this equality  when we remove it from the definition of eigenvalue in \cite[Definition~1.1]{IO}.) In the noncommutative setting, a flow is a C* dynamical system $(A,\R,\alpha)$ consisting of a (unital) C*-algebra $A$ and a one-pa\-ra\-meter group of automorphisms $\alpha\colon \R\to\Aut(A)$. The action $\alpha$ is said to be {\em minimal} (or, equivalently, we say  that $A$ is $\alpha$-simple) if $A$ has no nontrivial invariant ideals. When $A=C(Y)$ is a commutative (unital) C*-algebra, a flow $\alpha$ on $A$ is induced by a flow $T$ on $Y$, that is, $\alpha_t(f)=f\circ T^t$ for all $f$ in $A$ and for all $t$ in $\R$. Then $\alpha$ is minimal if and only if $T$ is minimal in the classical sense.

The aim in this paper is to extend some of the results in \cite{IO} to the noncommutative case, in other words, given a minimal C*-dynamical system $(A,\R,\alpha)$, we try to relate the values of $t$ for which the crossed product induced by the time $t$ automorphism $\alpha_t$ is not simple, with the eigenvalues of (a restriction of) $\alpha$.  We are unable to answer this question in general. Our partial results are contained in Section~2 of this paper. It is natural to ask what is the corresponding problem for von Neumann algebras. It turns out that minimality in the C*-algebra setting corresponds to central ergodicity in the W*-algebra setting (recall that a W* dynamical system $(N,\R,\theta)$ is said to be {\em centrally ergodic} if the restriction of $\theta$ to the center of $N$ is ergodic). Indeed, if $(A,\R,\alpha)$ is a minimal C* dynamical system then the crossed product $A\rtimes\sb\alpha\R$ is minimal if and only if the strong Connes spectrum of $\alpha$ is $\R$, cf.\ \cite{K}; whilst if $(N,\R,\theta)$ is a centrally ergodic W* dynamical system then the crossed product $N\rtimes\sb\theta\R$ is a factor if and only if the Connes spectrum of $\theta$ is $\R$, cf. \cite[Corollay~XI.2.8, pg~336]{T2}. In Section~1 we discuss this version of the problem, which we are able to solve satisfactorily. We conclude each of the sections with some open problems.

     The author gratefully acknowledge support and encouragement from professors
T.\ Giordano, D.\ Handelman and V.\ Pestov. In particular, we are grateful to T.\ Giordano for helpful explanations on the theory of von Neumann algebras and for many stimulating conversations.

%%
%% SECTION W*
%%

\section{W* dynamical systems}

Let $(N,\R,\theta)$ be a W* dynamical system. We say that a real number $s$ is an eigenvalue for $\theta$ if there exists a nonzero $a$ in $N$ such that, for all $t$ in $\R$, we have $\theta_t(a)=e^{i st}a$. When $\theta$ is ergodic (that is, the set of fixed points of $\theta$ is $\C 1$), the set of eigenvalues, which we denote by $\E(\theta)$, is a subgroup of $\R$, cf.\ \cite{S}.

Our first lemma is known for the case when the underlying measurable space of the center of the von Neumann algebra in question is a probability space, see \cite[Lemma~12.1.1, pg 326]{CFS}.  The proof there, however, can not be adapted to a more general situation. Our proof here is essentially contained in the proof of \cite[Theorem~10.6]{T}.

\begin{lemma}\label{Sinai}
Let $(N,\R,\theta)$ be a centrally ergodic W* dynamical system and let  $T$ be a strictly positive real number. If $\theta_T$ is not centrally ergodic then there exists a nonzero eigenvalue $s$ of the restriction of $\theta$ to the center of $N$, such that $e^{isT}= 1$.
\end{lemma}
\begin{proof}
Suppose that $\theta_T$ is not centrally ergodic. Then ${\mathcal A}=Z(N)^{\theta_T}$ is a commutative von Newmann algebra which is not reduced to the scalars. Furthermore, $\mathcal A$ is $\theta$-invariant and the action of $\theta$ on $\calA$ is ergodic and periodic. Hence the action of $\theta$ on $\calA$ is transitive. Thus there exists $T_0>0$ such that the action of $\theta$ on $\calA$  is isomorphic to the canonical action of $\R$ on $L\sp\infty(\R\slash T_0\Z)$. Since $s=\frac{2\pi}{T_0}$ is an eigenvalue for this action (with eigenfunction defined by $f(t+T_0\Z)=e^{ist}$, for all $t$ in $\R$), then we have that $s$ is an eigenvalue of the action of $\theta$ on $\calA$, and so, $s$ is an eigenvalue of the action of $\theta$ on $Z(N)$. Since the action of $\theta$ on $\calA$ is periodic with period $T_0$, there  is $k$ in $\Z$ such that $T=kT_0$.  Hence
$s=\frac{2\pi}{T_0}=\frac{2\pi k}{T}$ is a nonzero eigenvalue for the restriction of $\theta$ to  $Z(N)$ and $e^{isT}=e^{i2\pi k}=1$, as wanted.
\end{proof}

The next result can be regarded as the W* version of \cite[Proposition~1.5]{IO}.

\begin{proposition}\label{W*}
Let $(N,\R,\theta)$ be a centrally ergodic W* dynamical system and let $\widetilde{\theta}$ denote the restriction of $\theta$ to the center of $N$. Consider the map with domain ${\ts{\frac{1}{2\pi}}} \Q \otimes \E\left( \ttheta \right)$  and codomain $\R$ defined by ${ \ts{\frac{r}{2\pi}}}\otimes s \mapsto {\ts{ \frac{rs}{2\pi} }}$. This map is a $\Q$-linear monomorphism with range equal to
\[
 \Bigl\{0\Bigr\} \cup
 \left\{ t\in\R\setminus\{0\}\colon \theta_{\frac{1}{t}} \text{ is not centrally ergodic}\right\}.
\]
Hence the set above is a  $\Q$-linear subspace of $\R$ isomorphic to ${\ts{\frac{1}{2\pi}}} \Q \otimes \E\left( \ttheta \right)$.
\end{proposition}
\begin{proof}
It is clear that the map $\frac{r}{2 \pi}\otimes s\mapsto \frac{rs}{2 \pi}$ is a $\Q$-linear monomorphism. We only need to show that the range is as proposed. Let $r=\frac{p}{q}$ be a nonzero rational number and let $s$ be a nonzero eigenvalue for $\ttheta$. Put $T=\frac{rs}{2 \pi}=\frac{ps}{2\pi q}$. We show that $(N,\Z,\theta\sb\frac{1}{T})$ is not centrally ergodic. Since $\E\left( \ttheta \right)$ is a group, we get that $ps$ is also an eigenvalue. Therefore, there exists a nonzero $a$ in the center of  $N$ such that $\theta_t(a)=e^{i   pst}a$, for all $t\in\R$. Notice that $a$ is not a scalar since $ps\not = 0$. Hence
\[
   \theta\sb\frac{1}{T}(a)=\theta\sb\frac{2\pi q}{ps}(a)=e^{i 2\pi q}a=a.
\]
Thus $(N,\Z,\theta\sb\frac{1}{T})$ is not centrally ergodic. Conversely, assume that $(N,\Z,\theta\sb\frac{1}{T})$ is not centrally ergodic for some $T\not =0$. By Lemma~\ref{Sinai}, there is $0 \not =s\in \E\left( \ttheta \right)$
such that $e^{is\frac{1}{T}}=1$.  Therefore $s\frac{1}{T}=2\pi k $ for some $k\in \Z$ and so ${T}=\frac{1}{2\pi k} s$,  as desired.
\end{proof}

We now prove the following easy lemma.

\begin{lemma}\label{innerW*}
An  inner automorphism on a von Neumann algebra $N$ is centrally ergodic if and only if $N$ is a factor.
\end{lemma}
\begin{proof}
Let $\theta$ be an inner automorphism on a von Neumann algebra $N$.
If $N$ is a factor then it is clear that $\theta$ is centrally ergodic. To prove the converse,  suppose that $\theta$ is centrally ergodic.  Let $a$ be an element in the center of $N$. Since $\theta$ is inner, it follows  that $\theta(a)=a$. As $\theta$ is centrally ergodic, we conclude that $a$ is a scalar multiple of the identity, as wanted.
\end{proof}

We are ready to prove the main result of this section.

\begin{theorem}\label{THM-W*}
Let $(N,\R,\theta)$ be a centrally ergodic W* dynamical system where $N$ is not a factor. Let $t$ be a nonzero real number.  Denote by $\ttheta$ the restriction of $\theta$ to the center of $N$. The following statements are equivalent.
\begin{enumerate}
  \item The W* dynamical system $(N,\Z,\theta_t)$ is not centrally ergodic.
  \item The von Neumann algebra $\cprod{N}{\theta_t}$ is not a factor.
  \item There exists $(r,s)$ in $\Q\times\E\left(\ttheta\right)$ such that $rst=2\pi$.
\end{enumerate}
\end{theorem}
\begin{proof}
  (1)$\Rightarrow$(2): This is well known, see e.g.\ \cite[Corollay~XI.2.8, pg~336]{T2} or \cite[Theorem~8.11.15, pg~362]{P}.

  (2)$\Rightarrow$(3): If $\cprod{N}{\theta_t}$ is not a factor then either $(N,\Z,\theta_t)$ is not centrally ergodic or $\CSpec(\theta_t)\not=\T$, cf.\  \cite[Theorem~8.11.15, pg~362]{P}. If $(N,\Z,\theta_t)$ is not centrally ergodic, we use Proposition~\ref{W*} to conclude that $\frac{1}{t}=\frac{r}{2\pi}s$ for some $(r,s)$ in $\Q\times\E\left(\ttheta\right)$. Hence $rst=2\pi$, and we are done. Else, assume that $(N,\Z,\theta_t)$ is centrally ergodic and $\CSpec(\theta_t)\not=\T$. Then $\CSpec(\theta_t)\sp\perp\not=\{0\}$. As $\T$ is compact, we may use \cite[Theorem~XI.2.9(ii), pg~336]{T2} to conclude that $\theta_{t}^{n}=\theta_{nt}$ is inner for some nonzero integer $n$ in $\CSpec(\theta_t)\sp\perp\subset\Z$. Since $N$ is not a factor, we get that  $\theta_{nt}$ is not centrally ergodic, cf.\ Lemma~\ref{innerW*}. Another application of Proposition~\ref{W*} completes the proof.

  (3)$\Rightarrow$(1): Suppose there exists  $(r,s)$ in $\Q\times\E\left(\ttheta\right)$ such that $rst={2\pi}$. Then
 $\frac{1}{t}=\frac{r}{2\pi}s$ and so, by Proposition~\ref{W*}, the W* dynamical system $(N,\Z,\theta_t)$ is not centrally ergodic.
\end{proof}

 Recall that a von Neumann algebra $M$ of type $\rm III$ induces a (unique, up to conjugacy) flow on a von Neumann algebra of type $\rm{II}\sb\infty$, which is called the {\em noncommutative flow of weights of $M$}  or the {\em associated covariant system of $M$}, cf.\ \cite[Definition~XII.1.3, pg~368]{T2}.
As an example, we specialize Theorem~\ref{THM-W*} to the type III$\sb\lambda$ case, $0<\lambda<1$, to get the following.

\begin{corollary}\label{W*Cor}
Let $M$ be a factor of type {\rm III}$\sb\lambda$, $0<\lambda<1$, with associated covariant system $(N,\R,\theta)$.  Let $t$ be a nonzero real number. The following statements are equivalent.
\begin{enumerate}
  \item The W* dynamical system $(N,\Z,\theta_t)$ is not centrally ergodic.
  \item The von Neumann algebra $\cprod{N}{\theta_t}$ is not a factor.
  \item There exists a rational number $r$ such that $t=r\log\lambda$.
\end{enumerate}
\end{corollary}
\begin{proof}
Let $\ttheta$ denote the restriction of $\theta$ to the center of $N$. In this case ${\ts{\frac{2\pi}{\log\lambda}}}\Z=T(M)=\E\left(\ttheta\right)$, see \cite[Theorem~XII.1.6, pg~369]{T2} or \cite[28.11, pg~425]{Str}. Therefore, using Proposition~\ref{W*}, if $t$ is a nonzero real number  then $(N,\Z,\theta_t)$ is not centrally ergodic if and only if there is $r\sp\prime$ in $\Q$ and $n$ in $\Z$ such that $\frac{1}{t}=\frac{r\sp\prime}{2\pi}\frac{2\pi n}{\log\lambda}$ if and only if $t=r\log\lambda$, where $r=\frac{1}{r\sp\prime n}$ is rational.
\end{proof}

Let $\left( N,\R,\theta\right)$  be a centrally ergodic W* dynamical system.
It could be the case that, for all nonzero $t$, the crossed product induced by the time $t$ map  $\theta_t$ is a  factor: for example, if $(N,\R,\theta)$ is the associated covariant system of a factor of type $\rm III\sb{0}$ (because $\E\left(\ttheta\right)=\bigl\{0\bigr\}$, see  \cite[Theorem~XII.1.6, pg~369]{T2} or \cite[28.11, pg~425]{Str}).  On the other hand, it could be the case that, for every real number $t$, the crossed product induced by the time $t$ map $\theta_t$ is not a factor: for example, if  $\rcprod{N}{\theta}$ is semisimple, where $N$ is a properly infinite semifinite von Neumann algebra which admits a faithful semifinite normal trace $\tau$ such that $\tau\circ\theta_t=e^{-t}\tau$, for all $t$ in $\R$ (because $\E\left(\ttheta\right)=\R$, cf.\ \cite[Theorem~8.6]{T}).

Theorem~\ref{THM-W*} is no longer  valid when $N$ is a factor. The author is grateful to Professor A.\ Kishimoto for communicating to us one of his unpublished examples. We adapt his example here to the W* setting. Let $N$ be the von Neumann algebra generated by two  unitary operators $u$ and $v$ satisfying $uv=e^{2\pi i s}vu$, where $s$ is an irrational number.  Let $\theta$ be the flow on $N$ defined by
$\theta_t(u)= e^{2\pi i t }u$ and $\theta_t(v)= e^{2\pi i s t}v$, for all $t$ in $\R$.
 Then $\theta$ is centrally ergodic and $\E\left(\ttheta\right)=\{0\}$ (because $N$ is a factor). However, the crossed product  $\cprod{N}{\theta_1}$ induced by the time $1$ automorphism $\theta_1$ is not a factor (because $\theta_1$ is inner). It is worth mentioning that, in fact, the von Neumann algebra  $\rcprod{N}{\theta}$ is a factor. Therefore, this example satisfies a stronger condition than the one required in Theorem~\ref{THM-W*}.

One may ask if a similar result to Proposition~\ref{W*} is true if we substitute the centrally ergodic condition by ergodicity.  The example above suggests we should substitute the invariant $\E\left(\ttheta\right)$ by $\E\left(\theta\right)$. We remark that, by results of S\o rmer \cite{S}, if $\theta$ is a flow on $N$ which has kernel different from $\{ 0\}$ (or, equivalently, $\mathrm{Sp}(\alpha)\not=\R$, cf.\ \cite[Theorem~3.2]{S}), then $N$ must be abelian \cite[Theorem~3.5]{S}. This case, of course, is covered in Proposition~\ref{W*}. We conclude this section with some open questions.

\smallskip
\noindent{\bf Problem 1:} Suppose that $(N,\R,\theta)$ is a centrally ergodic flow. If $N$ is a factor, characterize the values of $t$ for which the crossed product associated to the time $t$ automorphism $\theta_t$ is a factor.

\smallskip
\noindent {\bf Problem 2:} Suppose that $(N,\R,\theta)$ is an ergodic flow. Characterize the values of $t$ for which the time $t$ automorphism $\theta_t$ is ergodic.
\smallskip

%%
%% SECTION C*
%%
\section{C* dynamical systems}

Let $(A,\R,\alpha)$ be a unital C* dynamical system. We say that a real number $s$ is an eigenvalue for $\alpha$ if there exists a nonzero $a$ in $A$ such that, for all $t$ in $\R$, it follows that $\alpha_t(a)=e^{ist}a$. We denote by $\E(\alpha)$ the set of eigenvalues of $\alpha$. In this section, we only consider eigenvalues of flows on commutative C*-algebras. If $A=C(Y)$ is commutative and $\alpha$ is induced by  a flow $T$ on $Y$, one may check that the eigenvalues for $\alpha$ are the same as the eigenvalues of $T$ in the classical sense. (We remark that in \cite[Definition~1.1]{IO} and some places in the literature,  a $2\pi$ term appears in the classical definition of eigenvalue. We compensate for such constant in the results below.)

Using results by Olesen and Pedersen, we obtain the following.

\begin{proposition}\label{Prime}
Let $(A,\R,\alpha)$ be a C* dynamical system. Assume that for all
$s$ in $\R$ and for all nonzero ideal $I$ of $A$,
$I\cap\alpha_s(I)\not = \{ 0 \}$. If $\alpha$ is minimal then, for
all nonzero real number $t$, the automorphism $\alpha_t$ is
minimal.
\end{proposition}
\begin{proof}
Let $t$ be a nonzero real number.  If  $\alpha_t$ is not minimal
then  there exists a nontrivial $\alpha_t$-invariant  ideal $J$ of
$A$.  Hence $J$ is invariant under $G_0=t\Z$. Since $\R\slash G_0$
is compact, we may use   \cite[Proposition~2.2]{OP2} to conclude
that $J$ contains a nonzero $\alpha$-invariant ideal $I$, and so,
the ideal $I$ is nontrivial because it is contained in the
nontrivial ideal $J$. This contradicts the minimality of $\alpha$
and completes the proof.
\end{proof}

We remark that if  $A$ is prime then the hypothesis of the
proposition is satisfied. Also, this hypothesis is equivalent to
say that the dual action $\widehat{\alpha}$ has full Connes
spectrum, cf.\ \cite[Lemma~3.2]{OP1}.

We can write the following result in a remarkable similarity to Proposition~\ref{W*}.

\begin{proposition}\label{C*}
Let $(A,\R,\alpha)$ be a unital and minimal  C* dynamical system and let $\talpha$ denote the restriction of $\alpha$ to the center of $A$. Assume that  $A$ is either commutative or  prime. Consider the map with domain ${\ts{\frac{1}{2\pi}}}\Q\otimes\E\left(\talpha\right)$
and codomain $\R$
given by ${\ts  { \frac{r}{2\pi}  \otimes s \mapsto \frac{rs}{2\pi} }}$. This map is a $\Q$-linear monomorphism with range equal to
\[                    \Bigl\{0\Bigr\} \cup
 \left\{ t\in\R\setminus\{0\}\colon \alpha_{\frac{1}{t}} \text{ is not minimal}\right\}.
\]
Hence the set above is a $\Q$-linear subspace of $\R$ isomorphic to ${\ts{\frac{1}{2\pi}}}\Q\otimes\E\left(\talpha\right)$.
\end{proposition}
\begin{proof}
If $A$ is commutative, this is \cite[Proposition~1.5]{IO}. For the case that $A$ is prime, we have that $Z(A)=\C 1$ and so $\E\left(\talpha\right)=\{0\}$. This, together with Proposition~\ref{Prime}, complete the proof.
\end{proof}

We now prove a statement analogous to Lemma~\ref{innerW*}.

\begin{lemma}\label{innerC*}
An  inner automorphism on a C*-algebra $A$ is minimal if and only if $A$ is simple.
\end{lemma}
\begin{proof}
Let $\alpha$ be an inner automorphism on a C*-algebra $A$. If $A$ is simple then it is clear that $\alpha$ is minimal.  To prove the converse, suppose that $\alpha$ is minimal. Let $I$ be an ideal of  $A$.  Then $I$ is $\alpha$-invariant because $\alpha$ is inner. Thus $I$ is trivial because $\alpha$ is minimal. Hence $A$ is simple.
\end{proof}

The following is the main result of this section.

\begin{theorem}\label{THM-C*}
Let $(A,\R,\alpha)$ be a unital and minimal C* dynamical system where $A$ is not simple. Assume that $A$ is either commutative or prime. Let $t$ be a nonzero real number. Denote by $\talpha$ the restriction of $\alpha$ to the center of $A$. The following statements are equivalent.
\begin{enumerate}
  \item The C* dynamical system $(A,\Z,\alpha_t)$ is not minimal.
  \item The C*-algebra $\cprod{A}{\alpha_t}$ is not simple.
  \item There exists $(r,s)$ in $\Q\times\E\left(\talpha\right)$ such that $rst=2\pi$.
\end{enumerate}
\end{theorem}
\begin{proof}
(1)$\Rightarrow$(2): This is well known, see eg.\ \cite[Theorem~6.5]{OP1} or \cite[Theorem~3.5 and Proposition~3.8]{K}.

(2)$\Rightarrow$(3): If $\cprod{A}{\alpha_t}$ is not simple then either $\alpha_t$ is not minimal or $\CSpec(\alpha_t)\not = \T$, cf.\  \cite[Theorem~6.5]{OP1}. If $\alpha$ is not minimal then, using Proposition~\ref{C*}, we can find a rational number $r$ and an eigenvalue $s$ in $\E\left(\talpha\right)$ such that $\frac{1}{t}=\frac{r}{2\pi}s$. Hence $rst=2\pi$, as desired. Else, assume that $\alpha_t$ is minimal and $\CSpec(\alpha_t)\not=\T$. Then $\CSpec\left(\alpha_t\right)\sp\perp\not=\left\{0\right\}$. As $\T$ is compact, we may use \cite[Theorem~4.5]{OP3} to find a nonzero $n$ in $\CSpec\left(\alpha_t\right)\sp\perp\subset\Z$ such that $\alpha_{t}^{n}=\alpha_{nt}$ is inner. Since  $A$ is not simple, we get that $\alpha_{nt}$ is not minimal, cf.\ Lemma~\ref{innerC*}. Another application of Proposition~\ref{C*} completes the proof.

(3)$\Rightarrow$(1): Suppose there exists $(r,s)$ in $\Q\times\E\left(\talpha\right)$ such that $rst=2\pi$. Then $\frac{1}{t}=\frac{r}{2\pi}s$ and so, by Proposition~\ref{C*}, $A$ is not $\alpha_t$-simple.
\end{proof}

As an example, we specialize to the case when $A$ is prime to get the following.

\begin{corollary}\label{primeNOsimple}
 Let $(A,\R,\alpha)$ be a unital C*-dynamical system. Assume that $A$ is prime but not simple. If $\alpha$ is minimal then $\cprod{A}{\alpha_t}$ is simple for all $0\not=t\in\R$.
\end{corollary}
\begin{proof}
 Follows from Theorem~\ref{THM-C*} since $Z(A)=\C1$ and so $\E(\talpha)=\{0\}$.
 \end{proof}

Theorem~\ref{THM-C*} (and Corollary~\ref{primeNOsimple}) fails for simple C*-algebras, as one can see from the (C*-version of) the example of A.\ Kishimoto described after Corollary~\ref{W*Cor}. We conclude this section with some open problems.

\smallskip
\noindent{\bf Problem 3:} Is Proposition~\ref{C*}  true if we erase the condition on $A$ being either commutative or prime?
\smallskip

\smallskip

\noindent{\bf Problem 4:} Suppose that $(A,\R,\alpha)$ is a minimal C* dynamical system and assume that $A$ is simple. Characterize the values of $t$ for which the crossed product associated to the time $t$ automorphism $\alpha_t$ is a simple C*-algebra.

%%%%%%%%%%%%%%%%%

\flushleft{\sc { \mbox{} \\
Department of Mathematics and Statistics.
University of Ottawa.  585 King Edward Ave. Ottawa, Ontario, K1N-6N5, Canada.\\
\smallskip
{\it E-mail address:} {\tt bitzaort@uottawa.ca}\\
\medskip
{\it Current Address:}\\ \smallskip
Centro de Investigaci\'on en Matem\'aticas. Universidad Aut\'onoma del Estado de Hidalgo. Carretera Pachuca-Tulancingo Km 4.5.
Pachuca de Soto, Hidalgo, 42090, M\'exico.\\
\smallskip
{\it Current e-mail address:} {\tt itza@uaeh.edu.mx}
}}
\end{document}